\documentclass{nlaauth}
\usepackage{color}

\begin{document}
\NLA{1}{6}{00}{28}{00}

\runningheads{D.\ Hezari, V.\ Edalatpour and D.\ K.\ Salkuyeh }
{Preconditioned GSOR iterative method}

\title{Preconditioned GSOR iterative method for a class of complex symmetric system of  linear equations}

\author{Davod Hezari, Vahid Edalatpour and Davod Khojasteh Salkuyeh\corrauth}

\address{Faculty of Mathematical Sciences, University of Guilan,\\
 Rasht, Iran. \\
 d.hezari@gmail.com, vedalat.math@gmail.com, khojasteh@guilan.ac.ir}

\corraddr{D. K. Salkuyeh, Faculty of Mathematical Sciences, University of Guilan, Rasht, Iran.}

\footnotetext[2]{khojasteh@guilan.ac.ir}

%\cgsn{Publishing Arts Research Council}{98--1846389}

\received{}
\revised{}
\noaccepted{}

\begin{abstract}
In this paper, we present a preconditioned variant of the generalized successive overrelaxation (GSOR) iterative method for
solving  a broad class of complex symmetric linear systems. We study conditions under which  the spectral radius of the iteration matrix of the preconditioned GSOR method is smaller than that of the GSOR  method  and determine the optimal values of iteration parameters. Numerical experiments are given to
verify the validity of the presented theoretical results and the effectiveness of the preconditioned GSOR method.
\end{abstract}

\keywords{complex linear systems; symmetric positive definite;  GSOR  method; preconditioning}

\section{INTRODUCTION}
Consider the system of linear equations of the form
\begin{equation}\label{1.1}
Au\equiv (W+iT)u= b, \hspace{.3cm} u,b\in\mathbb{C}^{n},
\end{equation}
where $i=\sqrt{-1}$ and $W,T\in\mathbb{R}^{n\times n}$ are symmetric positive
semidefinite matrices with at least one of them, e.g., $W$, being positive
definite. We assume that $T\neq0$, which implies that $A$ is non-Hermitian. Such systems arise in many problems such as quantum mechanics \cite{Dijk}, diffuse optical tomography \cite{Arridge}, structural dynamics \cite{Feriani}, FFT-based solution of certain time-dependent PDEs \cite{Bertaccini} and molecular scattering \cite{Poirier}. For more examples and additional references, the reader is referred to \cite{Benzi-B}.

The Hermitian/skew-Hermitian (HS) splitting of the matrix $A$ is
given by
\begin{equation}\label{1.2}
A=H+S,
\end{equation}
where
\[
H=\frac{1}{2}(A+A^H)=W \hspace{.3cm} {\rm and} \hspace{.3cm} S=\frac{1}{2}(A-A^H)=iT,
\]
with $A^H$ being the conjugate transpose of $A$. Based on the HS
splitting (\ref{1.2}), the HSS iteration method \cite{Bai-G} can
be straightforwardly applied to solve (\ref{1.1}). Bai et
al.\cite{Bai-B} recently proposed the use of the following
modified Hermitian and skew-Hermitian splitting (MHSS) method
which is  more efficient than the HSS iteration method for
solving the complex symmetric linear system
(\ref{1.1}):
\bigskip

\noindent \textbf{The MHSS iteration method.} \verb"Given an initial guess" $u^{(0)} \in \mathbb{C}^{n}$ \verb"and positive" \verb"constant"  $\alpha$ \verb"for" $k=0, 1, 2 \ldots$ \verb"until" \{$u^{(k)}$\} \verb"converges, compute"
\begin{equation}\label{1.3}
\left \{\begin{array}{ll}
(\alpha I+W)u^{(k+\frac{1}{2})}=(\alpha I-iT)u^{(k)}+b,\\
(\alpha I+T)u^{(k+1)}=(\alpha I+iW)u^{(k+\frac{1}{2})}-ib,
\end{array}\right.
\end{equation}
\verb"where" $I$ \verb"is the identity matrix".

\bigskip
In \cite{Bai-B}, Bai and coworkers proved that  the MHSS
iterative method is convergent for any positive constant $\alpha$.
Obviously both of the matrices $\alpha I+W$ and $\alpha I+T$ are
symmetric positive definite. Therefore, the two sub-systems
involved in each step of the MHSS iteration can be solved
effectively  by using the Cholesky factorization of the matrices
$\alpha I+W$ and $\alpha I+T$. Moreover, to solve both of the
sub-systems  in the inexact variant of the MHSS method one can
use a preconditioned conjugate gradient method. This is different
from the HSS iteration method, in which a shifted skew-Hermitian
linear sub-system with coefficient matrix $\alpha I+iT$ needs to
be solved at the second half-step of  every iteration. More recently, Bai et al. \cite{Bai-B-C1} proposed a preconditioned variant of the MHSS (PMHSS) for solving a class of complex symmetric systems of linear equations. It is necessary to mention that a potential difficulty with the HSS and MHSS iteration methods is the need to use complex arithmetic. Moreover, Axelsson et al. \cite{Axelson1} have presented a comparison of iterative methods to solve the complex symmetric linear system of equations (\ref{1.1}).

Letting $u=x+iy$ and  $b=p+iq$ where  $x,y,p,q\in
\mathbb{R}^{n}$, the complex linear system (\ref{1.1}) can be
rewritten as 2-by-2 block real equivalent formulation
\begin{equation}\label{1.4}
\mathcal{A}
  \begin{bmatrix}
    x \\
    y
  \end{bmatrix}
=
  \begin{bmatrix}
    p \\
    q
  \end{bmatrix},
\end{equation}
where
\[
\mathcal{A}=  \begin{bmatrix}
    W & -T \\
    T &W
  \end{bmatrix}.
\]

This  linear system can be formally regarded as a special case of the generalized saddle point problem \cite{Benzi1,Benzi2}. Recently, more efficient preconditioners for the real formulation (\ref{1.4}) have been proposed \cite{Bai-B-C2,Benzi-B,Giov,Howle}.
Rather than solving the original complex linear system
(\ref{1.1}), Salkuyeh et al. \cite{Salkuyeh} solved the real
equivalent system (\ref{1.4}) by the generalized successive
overrelaxation (GSOR) iterative method. By some numerical experiments, they have shown that the
performance of the GSOR method
is much more better than the MHSS method.

In order to solve the linear system (\ref{1.4}) by the GSOR method, the matrix $\mathcal{A}$ is split as
\[
\mathcal{A}=
\begin{bmatrix}
    W & 0 \\
    0 & W
  \end{bmatrix}-
\begin{bmatrix}
    0 & 0 \\
    -T & 0
  \end{bmatrix}-
\begin{bmatrix}
    0 & T \\
    0 & 0
  \end{bmatrix}.
\]
So, for $0\neq\alpha\in\mathbb{R}$, the GSOR method can be constructed as
follows:
\begin{equation}\label{1.5}
\begin{bmatrix}
  x^{k+1} \\
  y^{k+1}
  \end{bmatrix}
  =\mathcal{G}_{\alpha}
\begin{bmatrix}
  x^{k} \\
  y^{k}
  \end{bmatrix}+C_{\alpha}
  \begin{bmatrix}
  p \\
  q
  \end{bmatrix},
\end{equation}
where
\begin{eqnarray*}
\mathcal{G}_{\alpha}&=&
  \begin{bmatrix}
  W & 0 \\
  \alpha T & W
  \end{bmatrix}^{-1}\Bigg ((1-\alpha)
  \begin{bmatrix}
  W & 0 \\
  0 & W
  \end{bmatrix}+\alpha
  \begin{bmatrix}
  0 & T \\
  0 & 0
  \end{bmatrix}\Bigg), \\
&=&
  \begin{bmatrix}
  I & 0 \\
  \alpha S & I
  \end{bmatrix}^{-1}
  \begin{bmatrix}
 (1-\alpha)I & \alpha S \\
  0 & (1-\alpha)
  \end{bmatrix},
\end{eqnarray*}
is the iteration matrix, wherein  $S=W^{-1}T$, and
\begin{eqnarray*}
C_{\alpha}=\alpha
  \begin{bmatrix}
 W & 0 \\
  \alpha T & W
  \end{bmatrix}^{-1}.
\end{eqnarray*}

It is easy to see that (\ref{1.5}) is equivalent to
\begin{equation}\label{1.6}
\left \{\begin{array}{ll} Wx^{(k+1)}=(1-\alpha)Wx^{(k)} + \alpha Ty^{(k)}+ \alpha p,\\
Wy^{(k+1)}=-\alpha Tx^{(k+1)} + (1-\alpha)Wy ^{(k)}+ \alpha q,
\end{array}\right.
\end{equation}
where $x^{(0)}$ and $y^{(0)}$ are the initial approximations for $x$
and  $y$, respectively. Salkuyeh et al. \cite{Salkuyeh} proved that,
under certain condition on $\alpha$, the GSOR method is
convergent and determined the optimal value of the iteration
parameter $\alpha$ and the corresponding optimal convergence factor.
In the GSOR method two sub-systems with the coefficient matrix $W$
should be solved which can be done by the Cholesky factorization
or inexactly by the CG algorithm. Moreover, the right-hand side
of the sub-systems are real. Therefore, the solution of the
system can be obtained by the real version of the algorithms.

By decreasing the spectral radius of the iteration matrix $\mathcal{G}_{\alpha}$, the convergence rate of the GSOR method
can been improved. For this purpose, an effective method is to
transform the linear system (\ref{1.4}) into the preconditioned form
\[
\mathcal{P}\mathcal{A}
  \begin{bmatrix}
 x \\
  y
  \end{bmatrix}
=\mathcal{P}
  \begin{bmatrix}
 p \\
  q
  \end{bmatrix},
\]
where $\mathcal{P}\in\mathbb{R}^{2n \times 2n}$ is nonsingular. The
basic GSOR iterative method corresponding to the preconditioned
system will be referred to as the preconditioned GSOR iterative
method (PGSOR).

In this paper, we are going to consider the preconditioned GSOR method  with the following preconditioner
\[
\mathcal{P_{\omega}}=
  \begin{bmatrix}
 \omega I & I \\
 -I & \omega I
  \end{bmatrix},
\]
where $ 0<\omega\in\mathbb{R}$.  We  study condition on $\omega$ under which the spectral
radius of the  iterative matrix of the PGSOR method is smaller than that of the GSOR method and determine the optimal values of $\alpha$ and $\omega$. Moreover, we propose the approximation values for $\alpha$ and $\omega$ such that the performance of the corresponding PGSOR method is close to that of the optimal parameters. Finally,
numerical examples are presented to verify the theoretical results and the effectiveness of the  PGSOR
method. It is noteworthy that the motivation of choosing this preconditioner stems from \cite{Benzi-B} in which Benzi et al. have presented some examples of preconditioners for the  real formulation (\ref{1.4}) of the system (\ref{1.1}).

Throughout this paper, for a square matrix $X$, $\sigma(X)$ and $\rho(X)$ denote for the spectrum and the spectral radius of $X$, respectively. For a vector $z\in \Bbb{C}^n$, $\|z\|_2$ denotes for the Euclidean norm of $z$.

The paper is organized as follows. In Section 2 we discuss the
convergence of the PGSOR iterative method. Section
3 is devoted to some numerical experiments  to confirm the
theoretical results given in Section 2. Finally, in Section 4,
some concluding remarks are given.

\section{THE PGSOR ITERATIVE METHOD}

In this section, we express the PGSOR iterative
method and its convergence properties. Before presenting this
method,  we review the established theorems in
\cite{Salkuyeh} that relate to the convergence
properties of the GSOR method.
%%%%%%%%%%%%%%%%%%%%%%%%%%%%
\bigskip

\noindent \textsl{Theorem 2.1}\\
Let $W,T\in \mathbb{R}^{n \times n}$ be symmetric positive
definite and symmetric, respectively. Then, the  GSOR method to
solve Eq. (\ref{1.4}) is convergent if and only if
\[
0<\alpha<\frac{2}{1+\rho(S)},
\]
where $S=W^{-1}T$.

\bigskip
%%%%%%%%%%%%%%%%%%%%%%%%%%%%%%%%%%%%
%\bigskip

\noindent \textsl{Theorem 2.2}\label{t2}\\
 Let $W,T\in \mathbb{R}^{n \times n}$ be symmetric positive
definite and symmetric, respectively. Then, the optimal value of
the parameter $\alpha$ for the GSOR iterative method
(\ref{1.5})
 is given by
 \begin{equation}\label{1.7}
\alpha^*=\frac{2}{1+\sqrt{1+\rho(S)^{2}}},
\end{equation}
and the corresponding optimal convergence factor of the method is
given by
\begin{equation}\label{1.8}
\rho(\mathcal{G}_{\alpha^*})=1-\alpha^*=1-\frac{2}{1+\sqrt{1+\rho(S)^{2}}},
\end{equation}
where $S=W^{-1}T$.

\bigskip
%%%%%%%%%%%%%%%%%%%%%%%%%%%%%%%%%%%%%%%%%%%%%%%%%%
%\bigskip

\noindent \textsl{Corollary 2.1}\\
Let $W,T\in \mathbb{R}^{n \times n}$ be symmetric positive
definite and symmetric positive semi-definite, respectively.
Then, the GSOR method is convergent if and only if
$$0<\alpha<\frac{2}{1+\mu_{\max}(S)},$$
where $\mu_{\max}(S)$ is largest eigenvalue of $S=W^{-1}T$.
Moreover, the optimal value of iteration parameter $\alpha$ and
corresponding optimal convergence factor can be computed as
following
\begin{equation}\label{1.9}
\alpha^*=\frac{2}{1+\sqrt{1+\mu_{\max}(S)^{2}}}~~~~{\rm and}~~~~\rho(\mathcal{G}_{\alpha^*})=1-\alpha^*=1-\frac{2}{1+\sqrt{1+\mu_{\max}(S)^{2}}}.
\end{equation}

\bigskip
%%%%%%%%%%%%%%%%%%%%%%%%%%%%%%%%%%%%
%\bigskip

\noindent \textsl{Remark 2.1}\\
When the matrices $W$ and $T$ are symmetric positive
definite and symmetric positive semi-definite, respectively, if $\mu_{\max}(S)=0$ ($S=W^{-1}T$), then according to Eq.(\ref{1.9}) we have $\alpha^{*}=1$ and $\rho(\mathcal{G}_{\alpha^*})=0$. This means that the method would have the highest speed of convergence. Hence, hereafter, we assume $\mu_{\max}(S)\neq 0$.

\bigskip
%%%%%%%%%%%%%%%%%%%%%%%%%%%%%%%%%%%%
In order to speed up the convergence rate of the GSOR method, it is known
that the spectral radius of the iterative matrix
$\mathcal{G_{\alpha}}$ must be decreased. According to Theorem
\ref{t2}, the spectral radius of the iteration  matrix
$\mathcal{G_{\alpha}}$ tends to zero as the spectral radius of $S$ approaches to zero. For decreasing the spectral radius of $S$, an effective method is to precondition the linear system (\ref{1.4}). We apply $\mathcal{P}_\omega$ to the system (\ref{1.4}) to obtain the preconditioned linear system
\begin{equation}\label{1.10}
\tilde{\mathcal{A}}_\omega
  \begin{bmatrix}
 x \\
 y
  \end{bmatrix}=
\begin{bmatrix}
 \tilde{p} \\
  \tilde{q}
\end{bmatrix},
\end{equation}
with
\[
\tilde{\mathcal{A}}_\omega=
\begin{bmatrix}
 \tilde{W}_{\omega} & -\tilde{T}_{\omega} \\
 \tilde{T}_{\omega} & \tilde{W}_{\omega}
\end{bmatrix}
=
\begin{bmatrix}
 \omega W+T & -(\omega T-W) \\
 \omega T-W & \omega W+T
\end{bmatrix}
\qquad {\rm and} \qquad
 \begin{bmatrix}
 \tilde{p} \\
 \tilde{q}
\end{bmatrix}=
 \begin{bmatrix}
 \omega p+q \\
 \omega q-p
\end{bmatrix}.
\]
%%%%%%%%%%%%%%%%%%%%%%%%%%%%%%%%%%%%
\bigskip

\noindent \textsl{Remark 2.2}\\
Since $W$ and $T$ are symmetric positive definite and symmetric
positive semi-definite, respectively, we can easily conclude that
$\tilde{W}_{\omega}=\omega W+T$ and $\tilde{T}_{\omega}=\omega
T-W$ are symmetric positive definite and symmetric, respectively.

%\bigskip
%%%%%%%%%%%%%%%%%%%%%%%%%%%%%%%%%%%%
By splitting  the coefficient matrix $\tilde{\mathcal{A}}_\omega$
as
\[
\tilde{\mathcal{A}}_\omega=
 \begin{bmatrix}
\tilde{W}_{\omega} & 0 \\
 0 & \tilde{W}_{\omega}
\end{bmatrix}-
 \begin{bmatrix}
0 & 0 \\
 -\tilde{T}_{\omega} & 0
\end{bmatrix}-
 \begin{bmatrix}
0 & \tilde{T}_{\omega} \\
 0 & 0
\end{bmatrix},
\]
the GSOR method for solving (\ref{1.10}) is given by
\begin{equation}\label{1.11}
\begin{bmatrix}
  x^{k+1} \\
  y^{k+1}
\end{bmatrix}
=\tilde{\mathcal{G}}_{\alpha}(\omega)
 \begin{bmatrix}
  x^{k} \\
  y^{k}
\end{bmatrix}
+ \tilde{C}_\alpha(\omega)
 \begin{bmatrix}
  \tilde{p} \\
 \tilde{q}
\end{bmatrix},
\end{equation}
where
\begin{eqnarray*}
\tilde{\mathcal{G}}_{\alpha}(\omega)=
 \begin{bmatrix}
  I & 0 \\
 \alpha \tilde{S}_{\omega} & I
\end{bmatrix}^{-1}
 \begin{bmatrix}
  (1-\alpha)I & \alpha \tilde{S}_{\omega} \\
 0 & (1-\alpha)I
\end{bmatrix},\\
\end{eqnarray*}
is iteration matrix, wherein
$\tilde{S}_{\omega}=\tilde{W}_{\omega}^{-1}\tilde{T}_{\omega}$, and
\begin{eqnarray*}
\tilde{C}_{\alpha}(\omega)=\alpha
 \begin{bmatrix}
  \tilde{W} & 0 \\
\alpha \tilde{T} & \tilde{W}
\end{bmatrix}^{-1}.
\end{eqnarray*}

It is easy to see that (\ref{1.11}) is equivalent to
\begin{equation}\label{1.12}
\left \{\begin{array}{ll} \tilde{W}_{\omega}x^{(k+1)}=(1-\alpha)\tilde{W}_{\omega}x^{(k)} + \alpha \tilde{T}_{\omega}y^{(k)}+ \alpha \tilde{p},\\
\tilde{W}_{\omega}y^{(k+1)}=-\alpha \tilde{T}_{\omega}x^{(k+1)} + (1-\alpha)\tilde{W}_{\omega}y ^{(k)}+ \alpha \tilde{q},
\end{array}\right.
\end{equation}
where $x^{(0)}$ and $y^{(0)}$ are the initial approximations for $x$ and  $y$, respectively. In the PGSOR method two sub-systems with coefficient matrix $\tilde{W}_{\omega}=\omega W+T$ should be solved which can be done by the Cholesky factorization
or inexactly by the CG algorithm. Often $\omega W+T$
is better conditioned than $W$ itself.
%%%%%%%%%%%%%%%%%%%%%%%%%%%%%%%%%%%%
\bigskip

\noindent \textsl{Remark 2.3}\\
The PGSOR method is equivalent to apply the GSOR method for real equivalent formation of the new complex  system that is obtained by multiplying the complex number $(\omega-i)$ through both sides of the complex system (\ref{1.1}).

\bigskip
%%%%%%%%%%%%%%%%%%%%%%%%%%%%%%%%%%%%
%\bigskip

\noindent \textsl{Lemma 2.1} (see\cite{Salkuyeh})\\
 Let $W,T \in \mathbb{R}^{n \times n} $ be
symmetric positive definite and symmetric, respectively. Then,
the eigenvalues of the matrix $S=W^{-1}T$  are all real. If $T$
be a symmetric positive semi-definite matrix, then, the
eigenvalues of $S$ are all nonnegative.

\bigskip
%%%%%%%%%%%%%%%%%%%%%%%%%%%%%%%%%%%%%%%%%%%%%%%%%%%%
%\bigskip

\noindent \textsl{Lemma 2.2}\\
Let $W,T\in \mathbb{R}^{n \times n}$ be symmetric positive
definite and symmetric positive semi-definite, respectively. Also let $\omega$ be a
positive constant, $\tilde{W}_{\omega}=\omega W+T$ and
$\tilde{T}_{\omega}=\omega T-W$. If $\lambda$  is an eigenvalue of
$\tilde{S}_{\omega}=\tilde{W}_{\omega}^{-1}\tilde{T}_{\omega}$,
then, there is an eigenvalue $\mu$  of $S=W^{-1}T$ that
satisfies
\begin{equation}\label{1.13}
  \lambda=\frac{\omega\mu-1}{\omega+\mu}.
\end{equation}
Moreover
\begin{equation}\label{1.14}
  \rho(\tilde{S}_\omega)=\max\left\{\frac{1-\omega\mu_{\min}}{\omega+\mu_{\min}},
\frac{\omega\mu_{\max}-1}{\omega+\mu_{\max}} \right\},
\end{equation}
where $\mu_{\min}$ and $\mu_{\max}$ are the
smallest and largest eigenvalues of $S$, respectively.

\textit{Proof}\\
Let $(\lambda,x)$ be an  eigenpair of  $\tilde{S}_{\omega}$. Then $\tilde{S}_{\omega}x = \lambda x$, which is equivalent to
\[
(\omega T-W)x = \lambda(\omega W+T)x.
\]
Multiplying  both sides of the above equation by $W^{-1}$, we obtain
\[
(\omega-\lambda)Sx = (\omega \lambda+1)x.
\]
Evidently, we have $\lambda\neq\omega$. Then, from latter equation we get
\[
Sx=\frac{\omega\lambda+1}{\omega-\lambda}x,
\]
and hence, there is an eigenvalue $\mu$  of $S=W^{-1}T$ such that
\[
\mu=\frac{\omega\lambda+1}{\omega-\lambda},
\]
and from this we can easily obtain (\ref{1.13}).

The second part is a consequence of the fact that
\[
h(\mu)=\frac{\omega\mu-1}{\omega+\mu},
 \]
 is an increasing function with respect to variable $\mu$. \qquad $\Box$

\bigskip
%%%%%%%%%%%%%%%%%%%%%%%%%%%%%%%%%%%%
Since $\tilde{W}_\omega$ and $\tilde{T}_\omega$ are, respectively,   symmetric positive definite and symmetric matrices, analogously to Theorems 2.1 and 2.2 we can prove that the PGSOR iteration  is convergent if and only if
$$0<\alpha<\frac{2}{1+\rho(\tilde{S}_\omega)},$$
and the optimal value of the
iteration parameter $\alpha$ for the PGSOR iterative
method (\ref{1.11}) is given by
 \begin{equation}\label{1.15}
\alpha^*=\frac{2}{1+\sqrt{1+\rho(\tilde{S}_\omega)^{2}}},
\end{equation}
and moreover, the corresponding optimal convergence factor of the method is given by
\begin{equation}\label{1.16}
\rho(\mathcal{G}_{\alpha^*}(\omega))=1-\alpha^*=1-\frac{2}{1+\sqrt{1+\rho(\tilde{S}_\omega)^{2}}},
\end{equation}
where  $\tilde{S}_\omega=(\omega W+T)^{-1}(\omega T-W)$.

Note that, in the above relations, for the value of $\rho(\tilde{S}_\omega)$ one may use Eq. (\ref{1.14}).
%%%%%%%%%%%%%%%%%%%%%%%%%%%%%%%%%%%%%%%%%%%%

The main objective of applying $\mathcal{P}_\omega$ is to decrease the spectral radius of the iteration  matrix. This work is done by decreasing the spectral radius of $\tilde{S}_\omega$ in comparison to that of  $S$, then, from (\ref{1.8}) and (\ref{1.16}), we will get  $\rho(\tilde{\mathcal{G}}_{\tilde{\alpha}^*}(\omega))<\rho(\mathcal{G}_{\alpha^*})$. The necessary and sufficient conditions to achieve these desired results are formally stated in the following lemma and theorem.

%%%%%%%%%%%%%%%%%%%%%%%%%%%%%%%%%%%%%%%%%%%%
\bigskip

\noindent \textsl{Lemma 2.3}\\
Let $W,T\in \mathbb{R}^{n \times n}$ be symmetric positive
definite and symmetric positive semi-definite, respectively. Also let $\omega$ be a
positive constant, $\tilde{W}_{\omega}=\omega W+T$ and
$\tilde{T}_{\omega}=\omega T-W$. Then, the spectral radius of
$\tilde{S}_{\omega}=\tilde{W}_{\omega}^{-1}\tilde{T}_{\omega}$ is
smaller than that of $S=W^{-1}T$ if and only if
\begin{equation}\label{1.17}
 \max\bigg\{0,\frac{1-\mu_{\min}\mu_{\max}}{\mu_{\min}+\mu_{\max}}\bigg\}<\omega
\end{equation}
where $\mu_{\min}$ and $\mu_{\max}$ are the smallest and largest
eigenvalues of $S$, respectively.

\textit{Proof}\\
Suppose  $\lambda$ is an arbitrary eigenvalue of $\tilde{S}$.
Using Lemma 2.2, there is an  eigenvalue $\mu$ of $S$ such
that $\lambda=\frac{\omega\mu-1}{\omega+\mu}$. Moreover, from Lemma 2.1, we saw
that the eigenvalues of $S$ are nonnegative.
Then, the spectral radius of $\tilde{S}_{\omega}$ is smaller than
that of $S$ if and only if, for every eigenvalue $\mu$ of $S$, the
following inequality holds
\begin{equation}\label{1.18}
-\mu_{\max}<\frac{\omega\mu-1}{\omega+\mu}<\mu_{\max}.
\end{equation}
 It is easy to see that the right inequality of
(\ref{1.18}) holds if and only if $0<\omega$ and the left
inequality of (\ref{1.18}) holds if and only if
\begin{equation}\label{1.19}
\frac{1-\mu\mu_{\max}}{\mu+\mu_{\max}}<\omega\quad \forall \mu.
\end{equation}
Define
\[
f(\mu)=\frac{1-\mu\mu_{\max}}{\mu+\mu_{\max}}.
\]
Since  $f(\mu)$ is a decreasing function, then, the
inequality (\ref{1.19}) holds if and only if
\[
\frac{1-\mu_{\min}\mu_{\max}}{\mu_{\min}+\mu_{\max}}<\omega,
\]
which completes the proof. \qquad $\Box$

\bigskip
%%%%%%%%%%%%%%%%%%%%%%%%%%%%%%%%%%%%
%\bigskip

\noindent \textsl{Theorem 2.3}\\
Let $W,T\in \mathbb{R}^{n \times n}$ be symmetric positive
definite and symmetric positive semi-definite, respectively, and
$\omega$ be a positive constant. Let
$\mathcal{G}_{\alpha^*}$ and
$\tilde{\mathcal{G}}_{\tilde{\alpha}^*}(\omega)$ be the iteration matrices
of the GSOR and PGSOR methods, respectively. Then,
$\rho(\tilde{\mathcal{G}}_{\tilde{\alpha}^*}(\omega))<\rho(\mathcal{G}_{\alpha^*})$
if and only if
\begin{equation}\label{1.20}
\max\bigg\{0,\frac{1-\mu_{\min}\mu_{\max}}{\mu_{\min}+\mu_{\max}}\bigg\}<\omega,
\end{equation}
where $\mu_{\min}$ and $\mu_{\max}$ are the smallest and
largest eigenvalues of $S=W^{-1}T$, respectively.

\textit{Proof}\\
From (\ref{1.8}) and (\ref{1.16}), we can find that
$\rho(\tilde{\mathcal{G}}_{\tilde{\alpha}^*}(\omega))<\rho(\mathcal{G}_{\alpha^*})$
if and only if $\rho(\tilde{S}_\omega)<\rho(S)$, where
$\tilde{S}_\omega=(\omega W+T)^{-1}(\omega T-W)$, and by
the Lemma 2.3, the latter inequality holds if and only if
$\omega$ satisfies (\ref{1.20}). \qquad $\Box$

%\bigskip
%%%%%%%%%%%%%%%%%%%%%%%%%%%%%%%%%%%%%%%%%%%%%%%%

In the previous theorem, for every $\omega$ satisfying (\ref{1.20}), we saw that the spectral radius of the PGSOR method is smaller than that of the GSOR  method. Now, in the next Theorem, we obtain the optimal value of the
parameter  $\omega$ which minimizes the spectral
radius of the iterative matrix of the PGSOR method,
i.e.,
\[
\rho(\tilde{\mathcal{G}}_{\alpha^*}(\omega^*))=\min_{\omega} \rho(\tilde{\mathcal{G}}_{\alpha^*}(\omega)).
\]
%%%%%%%%%%%%%%%%%%%%%%%%%%%%%%%%%%%%
%\bigskip

\noindent \textsl{Lemma 2.4}\\
Let $W,T\in \mathbb{R}^{n \times n}$ be symmetric positive
definite and symmetric positive semi-definite, respectively,
$\tilde{W}_{\omega}=\omega W+T$ and $\tilde{T}_{\omega}=\omega
T-W$. Moreover, let $\mu_{\min}$ and $\mu_{\max}$ be the smallest
and largest eigenvalues of $S=W^{-1}T$, respectively, and $\omega$
be a positive parameter. Then, the optimal value of the parameter
$\omega$ which minimizes the spectral radius
$\rho(\tilde{S}_{\omega})$ of the matrix
$\tilde{S}_{\omega}=\tilde{W}_{\omega}^{-1}\tilde{T}_{\omega}$ is
given by
\begin{equation}\label{1.21}
\omega^*=\frac{1-\mu_{\min}\mu_{\max}+\sqrt{(1+\mu_{\min}^2)(1+\mu_{\max}^2)}}{\mu_{\min}+\mu_{\max}}.
\end{equation}

\textit{Proof}\\
By using Lemma 2.1, it is known that the eigenvalues $\mu$ of
$S$ are nonnegative. Then, according to Lemma 2.2, we can
write
\[
\rho(\tilde{S}_{\omega})=\max_{\mu\in \sigma(S)}\frac{|\omega\mu-1|}{\omega+\mu}.
\]
Note that
\[
g(\mu)=\frac{\omega\mu-1}{\omega+\mu},
\]
is an increasing function with respect to variable $\mu$. Now, if
$\mu_{\max}\leq\frac{1}{\omega}$, then for all $\mu\in \sigma(S)$, we
have $\omega\mu-1\leq0$. Hence
\begin{equation}\label{1.22}
 h(\omega)\equiv\max_{\mu\in\mu(S)}\frac{|\omega\mu-1|}{\omega+\mu}
=\frac{1-\omega\mu_{\min}}{\omega+\mu_{\min}}~~~~{\rm
if}~~~~\mu_{\max}\leq\frac{1}{\omega}.
\end{equation}
If $\frac{1}{\omega}\leq\mu_{\max}$, then $0\leq\omega\mu_{\max}-1$. First, we suppose $\mu_{\min}\neq0$
and  consider the following two cases:
  \begin{description}
    \item[case I:] If $\mu_{\min}\leq\frac{1}{\omega}$, we have
  $\omega\mu_{\min}-1\leq0$, which implies that
  \begin{equation}\label{1.23}
  h(\omega)=\max\left\{\frac{1-\omega\mu_{\min}}{\omega+\mu_{\min}},
  \frac{\omega\mu_{\max}-1}{\omega+\mu_{\max}} \right\}~~~~{\rm
  if}~~~~\mu_{\min}\leq\frac{1}{\omega}\leq\mu_{\max}.
  \end{equation}
    \item[case II:] If $\frac{1}{\omega}\leq\mu_{\min}$, then
  $0\leq\omega\mu_{\min}-1$, which implies that
  \begin{equation}\label{1.24}
  h(\omega)=\frac{\omega\mu_{\max}-1}{\omega+\mu_{\max}}~~~~{\rm
  if}~~~~\frac{1}{\omega}\leq\mu_{\min}.
  \end{equation}
  \end{description}
From (\ref{1.22}), (\ref{1.23}) and (\ref{1.24}), we can obtain the
following result
\begin{equation}\label{1.25}
h(\omega)=\left\{\begin{array}{lll}
\displaystyle{\frac{1-\omega\mu_{\min}}{\omega+\mu_{\min}},\hspace{3.5cm}\rm{for}~~~\omega\leq\frac{1}{\mu_{\max}}},\\
\\
\displaystyle{\max\left\{\frac{1-\omega\mu_{\min}}{\omega+\mu_{\min}},
\frac{\omega\mu_{\max}-1}{\omega+\mu_{\max}} \right\},~~~\rm{for}~~~\frac{1}{\mu_{\max}}\leq\omega\leq\frac{1}{\mu_{\min}}},\\
\\
\displaystyle{\frac{\omega\mu_{\max}-1}{\omega+\mu_{\max}},\hspace{3.5cm}\rm{for}~~~\frac{1}{\mu_{\min}}\leq\omega}.
\end{array}\right.
\end{equation}
Define
\[
f_1(\omega)=\frac{1-\omega\mu_{\min}}{\omega+\mu_{\min}}~~~{\rm and}~~~f_2(\omega)=\frac{\omega\mu_{\max}-1}{\omega+\mu_{\max}}.
\]
The functions $f_1(\omega)$ and $f_2(\omega)$ are  decreasing and
increasing, respectively. From this fact, we can conclude that if
$\omega^*$ be the minimum point of $h(\omega)$, then it belongs to
interval $(\frac{1}{\mu_{\max}},\frac{1}{\mu_{\min}})$ and must satisfy
\[
\frac{1-\omega\mu_{\min}}{\omega+\mu_{\min}}=\frac{\omega\mu_{\max}-1}{\omega+\mu_{\max}}.
\]
 By simplifying the above equation, we get
\begin{equation}\label{1.26}
 \omega^2(\mu_{\min}+\mu_{\max})+2\omega(\mu_{\min}\mu_{\max}-1)-(\mu_{\min}+\mu_{\max})=0.
\end{equation}
The roots of Eq. (\ref{1.26}) are
\[
\omega_{\pm}=\frac{1-\mu_{\min}\mu_{\max}\pm\sqrt{(1+\mu_{\min}^2)(1+\mu_{\max}^2)}}{\mu_{\min}+\mu_{\max}}.
\]
It is easy to observe that $\omega_{-}<0$, then
\begin{equation}\label{1.27}
\omega^*={\rm arg}\min_{\omega} h(\omega)=\frac{1-\mu_{\min}\mu_{\max}+\sqrt{(1+\mu_{\min}^2)(1+\mu_{\max}^2)}}{\mu_{\min}+\mu_{\max}}.
\end{equation}
Now, suppose $\mu_{\min}=0$. Then, from (\ref{1.22}) and (\ref{1.23}) we can deduce
\begin{equation}\label{1.28}
h(\omega) =\left\{\begin{array}{ll}
\displaystyle{\frac{1}{\omega},\hspace{3.7cm}\rm{for}~~~\omega\leq\frac{1}{\mu_{\max}}},\\
\\
\displaystyle{\max\left\{\frac{1}{\omega},
\frac{\omega\mu_{\max}-1}{\omega+\mu_{\max}} \right\},~~~~~\rm{for}~~~\frac{1}{\mu_{\max}}\leq\omega}.\\
\end{array}\right.
\end{equation}
In this case, it is easy to see that $\omega^*$ is obtained by setting $\mu_{\min}=0$ in  (\ref{1.27}). \qquad $\Box$

\bigskip
%%%%%%%%%%%%%%%%%%%%%%%%%%%%%%%%%%%%%%%%%%%
%\bigskip

\noindent \textsl{Theorem 2.4}\\
 Let $W,T\in \mathbb{R}^{n \times n}$ be symmetric positive
definite and symmetric positive semi-definite, respectively. Then, the optimal values of
the  parameters $\alpha$ and $\omega$ for the PGSOR
iterative method (\ref{1.11}) is given by
\begin{equation}\label{1.29}
\omega^*=\frac{1-\mu_{\min}\mu_{\max}+\sqrt{(1+\mu_{\min}^2)(1+\mu_{\max}^2)}}{\mu_{\min}+\mu_{\max}}~~~~{\rm
and}~~~~\alpha^*=\frac{2}{1+\sqrt{1+\xi^{2}}},
\end{equation}
and the corresponding optimal convergence factor of the method is given by
\begin{equation}\label{1.30}
\rho(\tilde{\mathcal{G}}_{\alpha^*}(\omega^*))=1-\alpha^*=1-\frac{2}{1+\sqrt{1+\xi^{2}}},
\end{equation}
where
$$\xi\equiv\rho(\tilde{S}_{\omega^*})=\frac{1-\omega^*\mu_{\min}}{\omega^*+\mu_{\min}}~\left(= \frac{\omega^*\mu_{\max}-1}{\omega^*+\mu_{\max}} \right),$$
 wherein $\mu_{\min}$ and $\mu_{\max}$ are the smallest and largest
eigenvalues of $S=W^{-1}T$, respectively, and $\tilde{S}_{\omega^*}=({\omega^*} W+T)^{-1}({\omega^*} T-W)$.

\textit{Proof}\\
According to  (\ref{1.15}) and (\ref{1.16}), for every positive parameter
$\omega$, we have
\[
\alpha^*=\frac{2}{1+\sqrt{1+\rho(\tilde{S}_\omega)^{2}}}~~~{\rm
and}~~~\rho(\tilde{\mathcal{G}}_{\alpha^*}(\omega))=1-\alpha^*=1-\frac{2}{1+\sqrt{1+\rho(\tilde{S}_\omega)^{2}}},
\]
where  $\tilde{S}_\omega=(\omega W+T)^{-1}(\omega T-W)$. Noticing that $\rho(\tilde{\mathcal{G}}_{\alpha^*}(\omega))$ is
an increasing function with respect to $\rho(\tilde{S}_\omega)$, it is minimized when $\rho(\tilde{S}_\omega)$ is minimum.
From the Lemma 2.4, we know that $\rho(\tilde{S}_\omega)$ is
minimized by
\[
\omega^*=\frac{1-\mu_{\min}\mu_{\max}+\sqrt{(1+\mu_{\min}^2)(1+\mu_{\max}^2)}}{\mu_{\min}+\mu_{\max}}.
\]
In the proof of Lemma 2.4, for the above $\omega^*$, we have seen that
\[
\frac{1-\omega^*\mu_{\min}}{\omega^*+\mu_{\min}}= \frac{\omega^*\mu_{\max}-1}{\omega^*+\mu_{\max}}.
\]
Then, by using Lemma 2.2, and replacing $\rho(\tilde{S}_{\omega^*})$ by $\xi$, the proof is completed. \qquad $\Box$

%\bigskip
%%%%%%%%%%%%%%%%%%%%%%%%%%%%%%%%%%%%%%%%%%%

Now, we can obtain an upper bound for the spectral radius of $\tilde{\mathcal{G}}_{\alpha^*}(\omega^*)$ that is formally stated in the following Corollary.
%%%%%%%%%%%%%%%%%%%%%%%%%%%%%%%%%%%%%%%%%%%
\bigskip

\noindent \textsl{Corollary 2.2}\\
Let the conditions of Theorem 2.4 hold. Then
\begin{equation}\label{1.31}
\rho(\tilde{\mathcal{G}}_{\alpha^*}(\omega^*))<\frac{\sqrt{2}-1}{\sqrt{2}+1}\approx 0.172,
\end{equation}
and
\[
\alpha^*\in (\frac{2}{1+\sqrt{2}}, 1)\approx (0.828, 1).
\]

\textit{Proof}\\
According to Eqs. (\ref{1.29}) and (\ref{1.30}),  the proof will be completed if we show the following inequality
\begin{equation}\label{1.32}
\xi=\rho(\tilde{S}_{\omega^*})=\frac{1-\omega^*\mu_{\min}}{\omega^*+\mu_{\min}}<1.
\end{equation}
Substituting  $\omega^*$ defined in  Eq. (\ref{1.29}) in the above inequality and  simplifying, yields
\[
\mu_{\max}+\mu_{\min}^2\mu_{\max}-\mu_{\min}\sqrt{(1+\mu_{\min}^2)(1+\mu_{\max}^2)}<1+\sqrt{(1+\mu_{\min}^2)(1+\mu_{\max}^2)}+\mu_{\min}^2,
\]
and therefore
\begin{equation}\label{1.33}
(1+\mu_{\min}^2)(\mu_{\max}-1)<(1+\mu_{\min})\sqrt{(1+\mu_{\min}^2)(1+\mu_{\max}^2)}.
\end{equation}
Now, if $\mu_{\max}\leq 1$, then the inequality (\ref{1.33}) holds, and so  (\ref{1.32}) is true. On the other hand, if  $\mu_{\max}> 1$, from  Eq. (\ref{1.33}) we get
\[
(1+\mu_{\min}^2)(\mu_{\max}-1)^2<(1+\mu_{\min})^2(1+\mu_{\max}^2),
\]
which is always true, and hence (\ref{1.32}) holds. \qquad $\Box$

\bigskip
%%%%%%%%%%%%%%%%%%%%%%%%%%%%%%%%%%%%%%%
Note that the upper bound (\ref{1.31}) for the spectral radius of $\tilde{\mathcal{G}}_{\alpha^*}(\omega^*)$ is a constant independent of both data and size of the problem.

From the proof of Corollary 2.2, we can see that all of the eigenvalues of $\tilde{S}_{\omega^*}$ are clustered in the interval $(-1,1)$. Hence,  when the spectral radius of $S$ is large, from (\ref{1.8}) and (\ref{1.31}), we can find that the spectral radius of $\mathcal{G}_{\alpha^*}$ is close to $1$, whereas, the spectral radius $\tilde{\mathcal{G}}_{\alpha^*}(\omega^*)$ is smaller than $0.172$. This demonstrates the superiority of the PGSOR method over the GSOR method. But, when the spectral radius of $S$ is smaller than $1$, it is expected that the performance of PGSOR method is similar to that of the GSOR method.

Since it may turn out to be difficult to find the optimal values of the parameters $\alpha$ and $\omega$, we propose to use the values of $\omega$ and $\alpha$ such that the performance of the corresponding PGSOR method  is close to that of the optimal parameters. To do so,  we consider the following two cases:\\[0.2cm]
Case I: If $\mu_{\min}=0$ (the smallest eigenvalue of $S=W^{-1}T$ is zero, e.g., when $T$ is $spsd$) then, with assumption $\omega=1$, from (\ref{1.28}), we have $h(\omega)=\rho(\tilde{S}_{\omega})=1$ and this means that all of eigenvalues of $\tilde{S}_{\omega}$ are clustered in the interval $[-1,1]$ and on the other hand, from (\ref{1.15}) and (\ref{1.16}), we have $\alpha^*=\frac{2}{1+\sqrt{2}}\approx 0.828$ and $\rho(\tilde{\mathcal{G}}_{\alpha^*}(\omega))=\frac{\sqrt{2}-1}{\sqrt{2}+1}\approx 0.172$, respectively.\\[0.2cm]
\noindent Case II: If $\mu_{\min}\neq0$ then, with assumption $\omega=1$, from (\ref{1.25}), it is easy to see that $h(\omega)=\rho(\tilde{S}_{\omega})<1$, hence, all of eigenvalues of $\tilde{S}_{\omega}$ are clustered in the interval $(-1,1)$. So, from (\ref{1.15}), we have $\alpha^*\in (\frac{2}{1+\sqrt{2}}, 1)\approx (0.828, 1)$. When $\mu_{\min}$ is very small and close to zero or $\mu_{\max}$ (the largest eigenvalue of $S=W^{-1}T$) is rather large, from (\ref{1.25}), it is easy to find that $h(\omega)=\rho(\tilde{S}_{\omega})\approx 1$, therefore, from (\ref{1.15}), it can be expected that $\alpha^*\approx \frac{2}{1+\sqrt{2}}\approx 0.828$.\\
\indent Therefore, for a broad class of problems, we can consider $\omega=1$ and  $\alpha=0.828$.

It is noteworthy that, if there exist real numbers $\beta$ and $\delta$ such that both matrices $\hat{W}:= \beta W+ \delta T$ and
$\hat{T}:= \beta T - \delta W$ are symmetric positive semidefinite with at least one of them positive definite, we can
first multiply both sides of (\ref{1.1}) by the complex number $\beta-i\delta$ to get the equivalent system
\[
(\hat{W} + i\hat{T})x=\hat{b}\quad \rm{with} \quad  \hat{{\it b}}:=(\beta-i\delta){\it b},
\]
and then employ the PGSOR iteration method to the  equivalent real
system that is obtained from the above system.
%%%%%%%%%%%%%%%%%%%%%%%%%%%%%%%%%%%%%%%%%%%%%%%%%%
\section{NUMERICAL EXAMPLES}
In this section, we use three  test problems  from \cite{Bai-B}
and an example of \cite{Bertaccini} to illustrate the  effectiveness of the preconditioned GSOR
iteration method for solving the equivalent real system
(\ref{1.4}). We also compare the performance of the PGSOR method with  the  HSS, MHSS  and GSOR methods, in terms of  both iteration count (denoted by IT) and CPU time (in seconds, denoted by CPU). The HSS and MHSS iterations are
employed to solve the  complex system (\ref{1.1}) and the
two other methods  to solve the equivalent real system
(\ref{1.4}). The two half-steps comprising each iteration  of  the HSS method are computed by the Cholesky factorization and LU decomposition of the coefficient matrices.
 In each iteration of  the MHSS, GSOR and PGSOR iteration
methods, we use the Cholesky factorization of the coefficient
matrices to solve the sub-systems. The CPU times reported  are
the sum of the CPU times for the convergence of the method and
the CPU times for computing the Cholesky factorization and LU decomposition. It is necessary to mention that to solve symmetric positive definite system of linear equations we have used the sparse Cholesky factorization incorporated with the symmetric approximate minimum degree reordering. To do so, we have used the \verb"symamd" command of M{\small ATLAB} Version 7.

All the numerical experiments  were computed in double precision
using some  M{\small ATLAB} codes  on a 64-bit 1.73 GHz intel
Q740 core i7 processor and 4GB RAM running widows 7. We use a null vector as an initial guess and
the  stopping criterion
\[
\frac{\|b-Au^{(k)}\|_{2}}{\| b\|_{2}}<10^{-6},
\]
is always used where $u^{(k)}=x^{(k)}+iy^{(k)}.$

\bigskip

%------------------------- Example 1 ------------------------------
\noindent \textsl{Example 3.1} (see \cite{Bai-B})\\
 Consider the linear system of equations
\[
\left[\left(K+\frac{3-\sqrt{3}}{\tau}I\right)+i\left(K+\frac{3+\sqrt{3}}{\tau}I\right)\right]x=b,
\]
where $\tau$ is the time step-size and $K$ is the five-point
centered difference matrix approximating the negative Laplacian
operator $L=-\Delta$ with homogeneous Dirichlet boundary
conditions, on a uniform mesh in the unit square $[0, 1]\times[0,
1]$ with the mesh-size $h=1/(m+1)$. The matrix $K\in
\mathbb{R}^{n\times n}$ possesses the tensor-product form
$K=I\otimes V_{m}+V_{m}\otimes I$, with $V_{m}=h^{-2}{\rm
tridiag}(-1,2,-1)\in \mathbb{R}^{m\times m}$. Hence, $K$ is an
${n\times n}$ block-tridiagonal matrix, with $n=m^{2}$. We take
\[
W=K+\frac{3-\sqrt{3}}{\tau}I~~~{\rm and}~~~T=K+\frac{3-\sqrt{3}}{\tau}I,
\]
and the right-hand side vector $b$ with its $j$th entry $b_{j}$
being given by
$$b_{j}=\frac{(1-i)j}{\tau(j+1)^{2}},\hspace{.3cm} j=1,2,\ldots,n.$$
In our tests, we take $\tau=h$. Furthermore, we normalize
coefficient matrix and right-hand side by multiplying both by
$h^{2}$.

\bigskip

 %------------------------- Example 2 ------------------------------

\noindent \textsl{Example 3.2} (see \cite{Bai-B})\\
 Consider the linear system of equations
$$\left[(-\omega^{2}M+K )+i(\omega C_{V}+C_{H}) \right]x=b,$$
where $M$ and $K$ are the inertia and the stiffness
matrices, $C_{V}$ and $C_{H}$ are the viscous and the hysteretic
damping matrices, respectively, and $\omega$ is the driving
circular frequency. We take $C_{H}=\mu K$ with $\mu$ a damping
coefficient, $M=I$ , $C_{V}=10I$ , and $K$ the five-point centered
difference matrix approximating the negative Laplacian operator
with homogeneous Dirichlet boundary conditions, on a uniform mesh
in the unit square $[0, 1]\times[0, 1]$ with the mesh-size
$h=1/(m+1)$. The matrix $K\in \mathbb{R}^{n\times n}$
possesses the tensor-product form $K=I\otimes V_{m}+V_{m}\otimes
I$, with $V_{m}=h^{-2}{\rm tridiag}(-1,2,-1)\in
\mathbb{R}^{m\times m}$. Hence, $K$ is an ${n\times n}$
block-tridiagonal matrix, with $n=m^{2}$. In addition, we set
$\omega=\pi$, $\mu=0.02$, and the right-hand side vector $b$ to be
$b=(1 + i)A\textbf{1}$, with $\textbf{1}$ being the vector of all entries equal to
$1$. As before, we normalize the system by multiplying both sides
through by $h^{2}$.

\bigskip

%------------------------- Example 3 ------------------------------

\noindent \textsl{Example 3.3} (see \cite{Bai-B})\\
 Consider the linear system of equations $(W+iT)x=b$, with
$$T=I\otimes V+V\otimes I~~~~{\rm and}~~~~W=10(I\otimes V_{c}+V_{c}\otimes I)+9(e_{1}e_{m}^{T}+e_{m}e_{1}^{T})\otimes I,$$
where $V={\rm tridiag}(-1,2,-1)\in \mathbb{R}^{m\times m}$,
$V_{c}=V-e_{1}e_{m}^{T}-e_{m}e_{1}^{T}\in \mathbb{R}^{m\times m}$
and $e_{1}$  and $e_{m}$  are the first and last unit vectors in
$\mathbb{R}^{m}$, respectively.We take the right-hand side vector
$b$ to be $b=(1 + i)A\textbf{1}$, with $\textbf{1}$ being the vector of all entries
equal to $1$.

Here $T$ and $W$ correspond to the five-point centered difference
matrices approximating the negative Laplacian operator with
homogeneous Dirichlet boundary conditions and periodic boundary
conditions, respectively, on a uniform mesh in the unit square
$[0, 1]\times[0, 1]$ with the mesh-size $h=1/(m+1)$.

\bigskip

%------------------------- Example 4 ------------------------------

\noindent \textsl{Example 3.4} (see \cite{Bertaccini,LPHSS})\\
We consider the complex Helmholtz equation
\[
-\triangle u+\sigma_1 u + i \sigma_2 u =f,
\]
where $\sigma_1$ and $\sigma_2$ are real coefficient functions,  $u$ satisfies Dirichlet boundary conditions
in $D = [0,1] \times [0, 1]$ and $i=\sqrt{-1}$.
We discretize the problem with finite differences on a $m\times m$ grid with mesh size $h = 1/(m + 1)$.  This leads to a system of linear equations
\[
\left((K+\sigma_1 I)+i \sigma_2 I\right)x=b,
\]
where $K=I\otimes V_{m}+V_{m}\otimes I$ is the discretization of $-\triangle$  by means of centered differences, wherein $V_{m}=h^{-2}{\rm tridiag}(-1,2,-1)\in \mathbb{R}^{m\times m}$. The right-hand side vector $b$ is taken to
be $b=(1 + i)A\textbf{1}$, with $\textbf{1}$ being the vector of all entries equal to $1$.
Furthermore, before solving the system we normalize the coefficient matrix and the right-hand side vector by multiplying both by $h^{2}$. For the numerical tests we set $\sigma_1=\sigma_2=100$.

For all the  examples, the optimal values of the parameters $\alpha$ and $\omega$ (denoted by $\alpha^{*}$  and $\omega^*$) are listed in Table \ref{Table1} for different values of $m$, that the former is  used in the HSS, MHSS, GSOR and PGSOR iterative methods and the latter is only used in the PGSOR iterative method. The experimentally found optimal parameters $\alpha^{*}$ for the HSS and MHSS are the ones resulting in the least numbers of iterations for the two methods for each of the numerical examples and those are presented in \cite{Bai-B} (expect for Example 3.4). The value of $\alpha^{*}$ for the GSOR method is obtained from (\ref{1.9}) and the  values of  $\alpha^{*}$ and $\omega^{*}$ for the PGSOR method are obtained from (\ref{1.24}) in which the largest eigenvalue of matrix $S$ $(\mu_{\max}(S))$ has been  estimated by a few iterations of the power method, and  also if $T$ is symmetric positive definite then the smallest eigenvalue of matrix $S$ $(\mu_{\min}(S))$ can be estimated by a few iterations of the inverse power method and if $T$ is symmetric positive semi-definite then $\mu_{\min}(S)=0$.

From  Table \ref{Table1},  as we expected  for all the examples, the value of $\alpha^{*}$ in the  PGSOR method  belongs to the interval $(0.828, 1)$, independently of the data and  problem size (in fact, this confirms Corollary 2.2) and  decreases as the  mesh-size $h$ decreases.
But note that the rate of decrease for  $\alpha^{*}$ in the  PGSOR method decreases as the mesh-size $h$ decreases and for large values of $m$ the value of $\alpha^{*}$ for Examples 3.2 and 3.4 becomes approximately constant and is approximately equal to $0.895$ and $0.869$, respectively. Moreover, for Examples 3.1 and 3.3,  the optimal parameter $\omega^{*}$ decreases  as the mesh-size $h$ decreases, whereas for Examples 3.2 and 3.4 it increases. Note that, for all the  examples, the rate of change of $\omega^{*}$ decreases as the  mesh-size $h$ decreases.

%------------------TABLE  1  optimal parameter ---------------------

 \begin{table}\label{Table1}
 \centering
 \caption{ The optimal parameters $\alpha^{*}$ and $\omega^{*}$ for the HSS, MHSS, GSOR and PGSOR methods.\label{Table1}} \vspace{-.2cm}
\begin{tabular}{llllllll}\vspace{-0.2cm}\\ \hline \vspace{-0.3cm} \\  %\cline{3-5}
Example   & Method  &         &Grid\\\cline{4-8}  \vspace{-0.3cm}\\\vspace{0.3cm}

          &         &         &$16\times 16$  & $32\times 32$ & $64\times 64$ & $128\times 128$ & $256\times 256$ \vspace{-0.3cm}\\  \hline \vspace{-0.3cm} \\ \vspace{0.1cm}

No. 1    &  HSS    & $\alpha^*$ & $0.81$        & $0.55$        & $0.37$        & $0.28$          & $0.20$\\\vspace{0.1cm}

          &  MHSS   & $\alpha^*$ & $1.06$        & $0.75$        & $0.54$        & $0.40$         &$0.30$\\\vspace{0.1cm}
          &  GSOR   & $\alpha^*$ & $0.550$       & $0.495$       & $0.457$       & $0.432$        &$0.428$\\\vspace{0.0cm}
          &  PGSOR  & $\alpha^*$ & $0.990$       & $0.987$       & $0.986$       & $0.984$        &$0.983$\\\vspace{0.1cm}
          &         & $\omega^*$ & $0.657$       & $0.624$       & $0.602$       & $0.590$        &$0.583$\\\vspace{0.1cm}

No. 2    &  HSS    & $\alpha^*$ & $0.42$        & $0.23$        & $0.12$        & $0.07$         &$0.04$\\\vspace{0.1cm}

          &  MHSS   & $\alpha^*$ & $0.21$        & $0.08$        & $0.04$        & $0.02$         &$0.01$\\\vspace{0.1cm}
          &  GSOR   & $\alpha^*$ & $0.455$       & $0.455$       & $0.455$       & $0.455$        &$0.455$\\\vspace{0.0cm}
          &  PGSOR  & $\alpha^*$ & $0.898$       & $0.896$       & $0.895$       & $0.895$        &$0.895$\\\vspace{0.1cm}
          &         & $\omega^*$ & $1.309$       & $1.323$       & $1.328$       & $1.330$        & $1.330$\\\vspace{0.1cm}

No. 3    &  HSS    & $\alpha^*$ & $4.41$        & $2.71$        & $1.61$        & $0.93$         &$0.53$\\\vspace{0.1cm}
          &  MHSS   & $\alpha^*$ & $1.61$        & $1.01$        & $0.53$        & $0.26$         &$0.13$\\\vspace{0.1cm}
          &  GSOR   & $\alpha^*$ & $0.908$       & $0.776$       & $0.566$       & $0.353$        & $0.199$\\
          &  PGSOR  & $\alpha^*$ & $0.982$       & $0.956$       & $0.918$       & $0.885$          & $0.864$\\\vspace{0.1cm}
          &         & $\omega^*$ & $3.001$       & $1.980$       & $1.437$       & $1.181$        &$1.063$\\\vspace{0.1cm}

No. 4    &  HSS    & $\alpha^*$ & $1.44$        & $0.77$        & $0.40$        & $0.21$         &$0.11$\\\vspace{0.1cm}
          &  MHSS   & $\alpha^*$ & $0.37$        & $0.09$        & $0.021$        & $0.005$         &$0.002$\\\vspace{0.1cm}
          &  GSOR   & $\alpha^*$ & $0.862$       & $0.862$       & $0.862$       & $0.862$        & $0.862$\\
          &  PGSOR  & $\alpha^*$ & $0.973$       & $0.970$       & $0.969$       & $0.969$          & $0.969$\\
          &         & $\omega^*$ & $2.587$       & $2.711$       & $2.745$       & $2.755$        &$2.757$\\\hline

\end{tabular}
\end{table}
%%%%%%%%%%%%%%%%%%%%%%%%%%%%%%%%%%%%%%%%%%%%%%%%%%%%%%%%%%%%%%%%%%
%%%%%%%%%%%%%%%%%%%%%%%%%%%%%%%%%%%%%%%%%%%%%%%%%%%%%%%%%%%%%%%

In Tables \ref{Table2}-\ref{Table5}, we have  reported numerical results for Examples 3.1-3.4. These tables present IT and CPU for the HSS, MHSS, GSOR and PGSOR methods. As seen, the PGSOR method outperforms GSOR and behaves much better than MHSS and HSS, especially when problem size increases.  The cause of such  performance is easily predictable when we observe the spectral radius of the iteration matrices of the GSOR and PGSOR  methods for the four examples (Table \ref{Table6}) that is obtained from (\ref{1.8}) and (\ref{1.30}), respectively. We see that the spectral radius of the iteration matrix of the PGSOR method is smaller than that of the GSOR method when the optimal parameters are employed.
For the numerical results of Example 3.3 reported in Table \ref{Table6}, we see that the spectral radius of the iteration matrix of the GSOR method grows rapidly with problem size, while that of the PGSOR method grows very slowly and is very smaller than that of the GSOR method when the
problem size is large.

In Fig. 1 we have compared the number of iterations for the convergence of the PGSOR method in conjunction with optimal values of $(\alpha,\omega)$  (denoted by $\textrm{PGSOR}_{\textrm{op}}$) and with $(\alpha,\omega)=(0.828,1)$ (denoted by $\textrm{PGSOR}_{\textrm{ap}}$) for all the examples. As seen, $(\alpha,\omega)=(0.828,1)$ can be considered as a reasonable approximation of the optimal value of $(\alpha,\omega)$.

%------------------TABLE  2 , NUMERICAL RESULTS FOR EXAMPLE 1--------------------

\begin{table}
\centering
\caption{ Numerical results for Example 3.1. \label{Table2}}
\vspace{-.2cm}
\begin{tabular}{llllllll}\vspace{-0.2cm}\\ \hline \vspace{-0.3cm} \\  %\cline{3-5}

Method    &  $m\times m$ & $16\times 16$ & $32\times
32$&$64\times 64$ & $128\times 128$ & $256\times
256$\vspace{-0.0cm}\\  \hline \vspace{-0.3cm} \\ \vspace{0.0cm}

HSS               & IT  & $44$   & $65$    &  $97$    & $136$   &  $191$\\\vspace{0.2cm}

                  & CPU & $0.07$ & $0.31$  &  $2.28$  & $18.97$ & $187.79$\\\vspace{0.0cm}

MHSS              & IT  & $40$   & $54$    &   $73$   &  $98$   &  $133$\\\vspace{0.2cm}

                  & CPU & $0.08$ & $0.30$  &   $1.53$ &  $8.50$ &  $52.57$\\\vspace{0.0cm}

GSOR              & IT  & $19$   & $22$    &   $24$   &  $26$   &  $27$\\\vspace{0.2cm}

                  & CPU & $0.04$ & $0.05$  &  $0.15$  &  $0.64$ & $2.88$\\\vspace{0.0cm}

PGSOR             & IT  & $4$    & $4$     &   $5$    &   $5$   & $5$\\
                  & CPU & $0.03$ & $0.04$  &  $0.07$  &  $0.23$ & $1.01$\\\hline

\end{tabular}
\end{table}
%%%%%%%%%%%%%%%%%%%%%%%%%%%%%%%%%%%%%%%%%%%%%%

%-------------------- TABLE  3 , NUMERICAL RESULTS FOR EXAMPLE 2 --------------------

\begin{table}
\centering
\caption{ Numerical results for Example 3.2. \label{Table3}}
\vspace{-.2cm}
\begin{tabular}{llllllll}\vspace{-0.2cm}\\ \hline \vspace{-0.3cm} \\  %\cline{3-5}

Method    &  $m\times m$ & $16\times 16$  & $32\times 32$  &$64\times 64$ & $128\times 128$ & $256\times 256$\vspace{-0.0cm}\\  \hline \vspace{-0.3cm} \\

HSS               & IT  & $86$   & $153$   &  $284$   & $540$   &  $1084$\\\vspace{0.2cm}

                  & CPU & $0.11$ & $1.15$  &  $5.77$  & $61.86$ & $692.11$\\\vspace{0.0cm}

MHSS              & IT  & $34$   & $38$    & $50$     & $81$    &$139$\\\vspace{0.2cm}

                  & CPU & $0.08$ & $0.22$  & $1.08$   & $7.24$ & $54.88$\\\vspace{0.0cm}

GSOR              & IT  & $26$   & $24$    &  $24$    & $23$    & $23$\\\vspace{0.2cm}

                  & CPU & $0.04$ & $0.06$  & $0.16$   & $0.60$  & $2.54$\\\vspace{0.0cm}

PGSOR             & IT  & $8$    & $7$     &   $8$    &   $8$   & $8$\\
                  & CPU & $0.03$ & $0.04$  &  $0.08$  &  $0.29$ & $1.26$\\\hline

\end{tabular}
\end{table}
%%%%%%%%%%%%%%%%%%%%%%%%%%%%%%%%%%%%%%%%%%%%%%

%--------------------TABLE  4 , NUMERICAL RESULTS FOR EXAMPLE 3 --------------------

\begin{table}
\centering
\caption{Numerical results for Example 3.3. \label{Table4}}
\vspace{-.2cm}
\begin{tabular}{llllllll}\vspace{-0.2cm}\\ \hline \vspace{-0.3cm} \\  %\cline{3-5}

Method    &  $m\times m$ & $16\times 16$  & $32\times 32$ &$64\times 64$ & $128\times 128$ & $256\times 256$\vspace{-0.0cm}\\  \hline \vspace{-0.3cm} \\

HSS               & IT  & $84$   & $137$   &  $223$   & $390$   &  $746$\\\vspace{0.2cm}

                  & CPU & $0.11$ & $0.60$  &  $4.77$  & $47.10$ & $556.16$\\\vspace{0.0cm}

MHSS              & IT  & $53$   & $76$    &  $130$   & $246$   & $468$\\\vspace{0.2cm}

                  & CPU & $0.11$ & $0.43$  & $2.84$   & $22.34$ & $194.36$\\\vspace{0.0cm}

GSOR              & IT  & $7$    & $11$    &$20$      & $35$    & $71$\\\vspace{0.2cm}

                  & CPU & $0.03$ & $0.05$  & $0.17$   & $1.05$  & $8.69$\\\vspace{0.0cm}

PGSOR             & IT  & $5$    & $6$     &   $7$    &   $8$   & $8$\\
                  & CPU & $0.03$ & $0.04$  &  $0.09$  &  $0.38$ & $1.67$\\\hline

\end{tabular}
\end{table}
%%%%%%%%%%%%%%%%%%%%%%%%%%%%%%%%%%%%%%%%%%%%%%
%------------------TABLE  5 , NUMERICAL RESULTS FOR EXAMPLE 4--------------------

\begin{table}
\centering
\caption{Numerical results for Example 3.4. \label{Table5}}
\vspace{-.2cm}
\begin{tabular}{llllllll}\vspace{-0.2cm}\\ \hline \vspace{-0.3cm} \\  %\cline{3-5}

Method    &  $m\times m$ & $16\times 16$ & $32\times
32$&$64\times 64$ & $128\times 128$ & $256\times
256$\vspace{-0.0cm}\\  \hline \vspace{-0.3cm} \\ \vspace{0.0cm}

HSS               & IT  & $25$   & $46$    &  $86$    & $161$   &  $300$\\\vspace{0.2cm}

                  & CPU & $0.04$ & $0.12$  &  $0.68$  & $5.04$ & $41.19$\\\vspace{0.0cm}

MHSS              & IT  & $30$   & $36$    &   $39$   &  $40$   &  $41$\\\vspace{0.2cm}

                  & CPU & $0.06$ & $0.16$  &   $0.59$ &  $2.41$ &  $10.94$\\\vspace{0.0cm}

GSOR              & IT  & $8$   & $8$    &   $8$   &  $8$   &  $7$\\\vspace{0.2cm}

                  & CPU & $0.03$ & $0.04$  &  $0.09$  &  $0.30$ & $1.25$\\\vspace{0.0cm}

PGSOR             & IT  & $5$    & $5$     &   $5$    &   $5$   & $5$\\
                  & CPU & $0.03$ & $0.04$  &  $0.07$  &  $0.24$ & $1.02$\\\hline

\end{tabular}
\end{table}
%%%%%%%%%%%%%%%%%%%%%%%%%%%%%%%%%%%%%%%%%%%%%%

%------------%%%%%%%%%------------Fig 1,2--------%%%%%%%%%%%%%%%%------------
\begin{figure}[!hbp]
\centerline{\includegraphics[height=5cm,width=6cm]{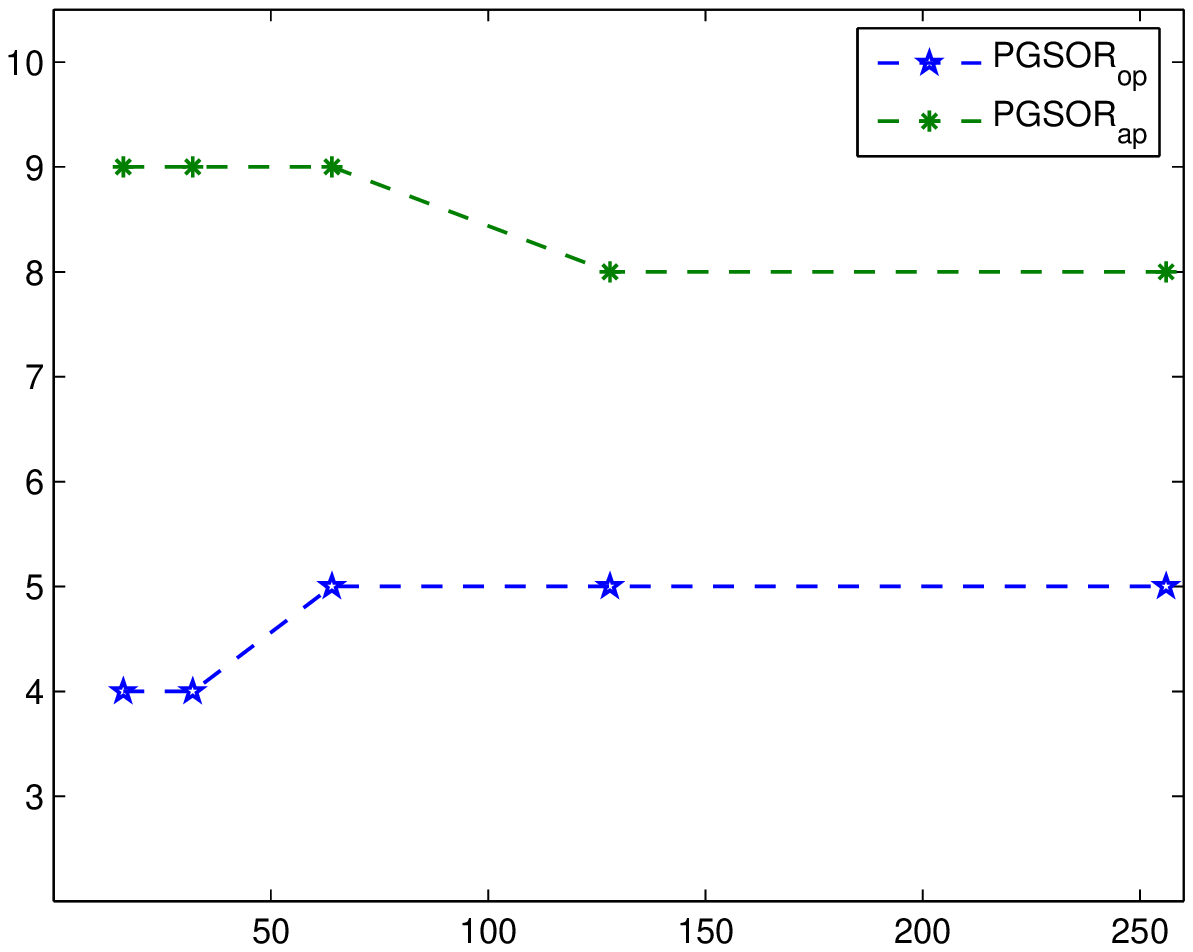}\includegraphics[height=5cm,width=6cm]{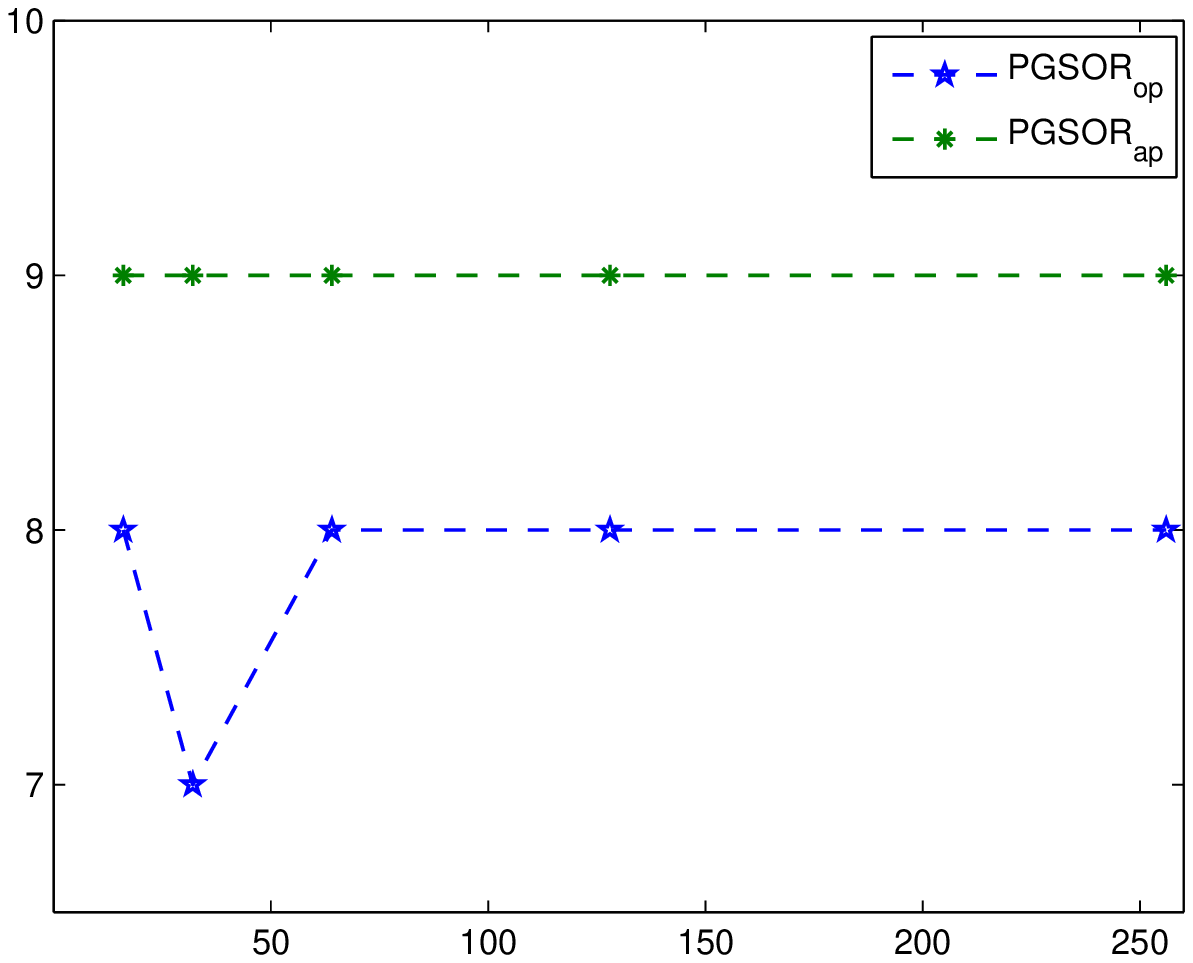}}
\centerline{\includegraphics[height=5cm,width=6cm]{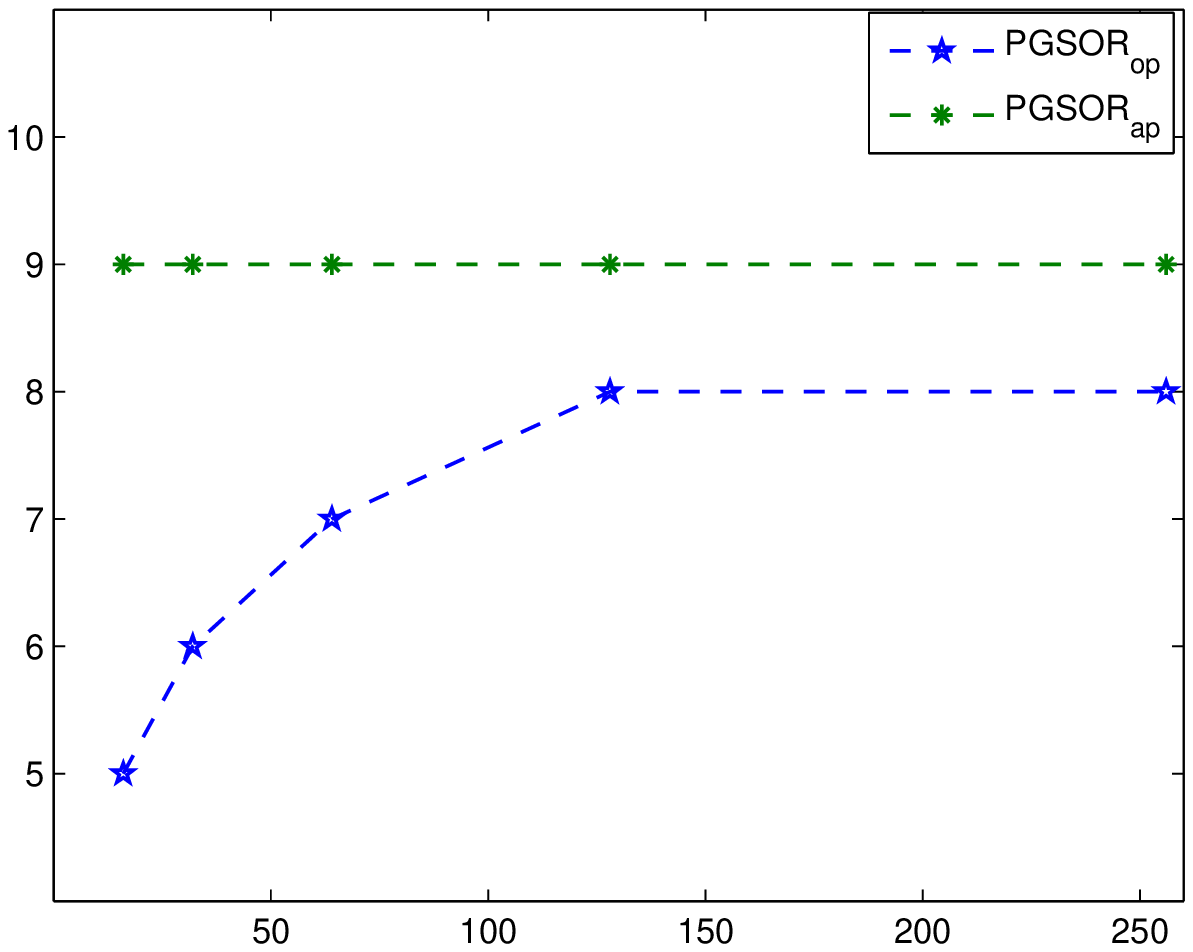}\includegraphics[height=5cm,width=6cm]{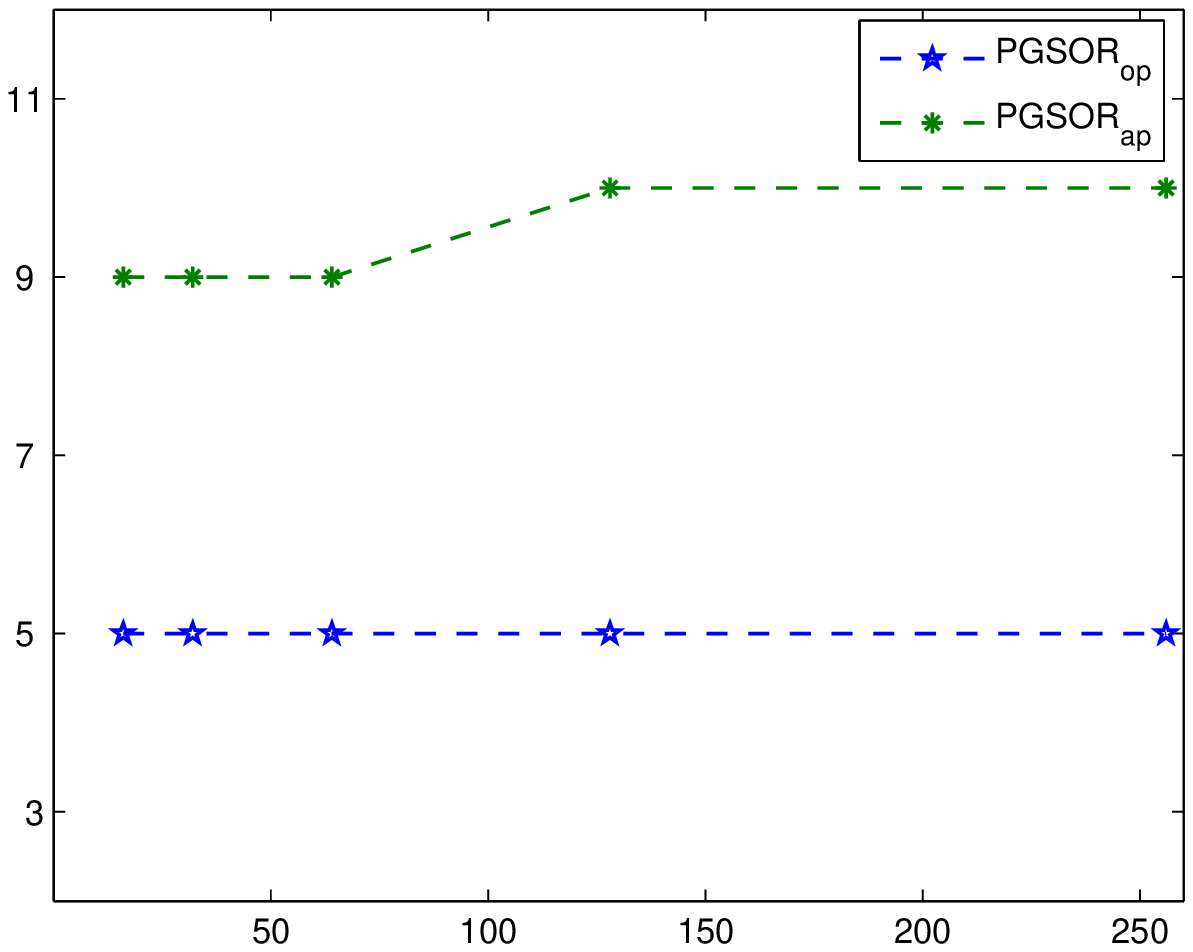}}
{\caption{Demonstration of IT versus $m$ for the PGSOR method with optimal value of the parameters $(\alpha,\omega)$ in Table \ref{Table1} and approximation parameters $\alpha=0.828$ and $\omega=1$; \textit{top-left}: Example 3.1, \textit{top-right}: Example 3.2,
\textit{down-left}: Example 3.3, \textit{down-right}: Example 3.4. }}\label{Fig2}\vspace{0cm}
\end{figure}

%\begin{figure}[!hbp]
%\centerline{\includegraphics[height=5cm,width=6cm]{Ex3.eps}\includegraphics[height=5cm,width=6cm]{Ex4.eps}}
%{\caption{Demonstration of IT versus $m$ for the PGSOR method with optimal values of the  parameters $(\alpha,\omega)$ in Table \ref{Table1} and  %$\alpha=0.828$ and $\omega=1$; \textit{left}: Example \ref{ex3}, \textit{right}: Example \ref{ex4}. }}\label{Fig3}\vspace{0cm}
%\end{figure}
%-----------------------------------------------------------------------------

%--------------------TABLE  5 , SPECTRAL RADIUS OF GSOR AND GSOR  --------------------

 \begin{table}
 \centering
 \caption{ Spectral radius of the iterative matrices of GSOR and PGSOR \label{Table6}}
\vspace{-.2cm}
\begin{tabular}{llllllll}\vspace{-0.2cm}\\ \hline \vspace{-0.3cm} \\  %\cline{3-5}
Example   &    & Grid\\\cline{3-7}  \vspace{-0.3cm}
\\\vspace{0.3cm}

          &         & $16\times 16$  & $32\times 32$ & $64\times 64$ & $128\times 128$ & $256\times 256$ \vspace{-0.3cm}\\  \hline \vspace{-0.3cm} \\ \vspace{0.2cm}

No. 1    &  $\rho(\mathcal{G}_{\alpha^*})$   & $0.450$        & $0.505$        &
           $0.543$        & $0.568$          & $0.572$\\\vspace{0.2cm}

          &  $\rho(\tilde{\mathcal{G}}_{\alpha^*}(\omega^*))$   & $0.010$       & $0.013$       & $0.014$       & $0.016$         & $0.017$\\\vspace{0.2cm}

No. 2    &  $\rho(\mathcal{G}_{\alpha^*})$   & $0.545$        & $0.545$        &
$0.545$        & $0.545$          & $0.545$\\\vspace{0.2cm}

          &  $\rho(\tilde{\mathcal{G}}_{\alpha^*}(\omega^*))$   & $0.102$       & $0.104$       & $0.105$       & $0.105$         & $0.105$\\\vspace{0.2cm}

No. 3    &  $\rho(\mathcal{G}_{\alpha^*})$   & $0.092$        & $0.224$        & $0.434$        & $0.647$          & $0.801$\\\vspace{0.2cm}
          &  $\rho(\tilde{\mathcal{G}}_{\alpha^*}(\omega^*))$   & $0.018$       & $0.044$       & $0.082$       & $0.115$         & $0.136$\\\vspace{0.2cm}

 No. 4    &  $\rho(\mathcal{G}_{\alpha^*})$   & $0.138$        & $0.138$        & $0.138$        & $0.138$          & $0.138$\\
          &  $\rho(\tilde{\mathcal{G}}_{\alpha^*}(\omega^*))$   & $0.027$       & $0.030$       & $0.031$       & $0.031$         & $0.031$\\\hline

\end{tabular}
\end{table}
%%%%%%%%%%%%%%%%%%%%%%%%%%%%%%%%%%%%%%%%%%%%%%
\section{CONCLUDING REMARKS}

In this paper we have presented a preconditioned variant of the
generalized successive overrelaxation (GSOR) iterative method to solve the equivalent
real formulation of complex linear system (\ref{1.1}), where $W$ is
symmetric positive definite and $T$ is symmetric positive
semi-definite. Convergence properties of the method have been also
investigated. Some numerical have been presented to show
the effectiveness of the method. Our numerical examples show that
our method is quite suitable for such problems. Moreover, the
presented numerical experiments show that the PGSOR method is
superior to GSOR, MHSS and HSS in terms of the iterations and CPU times.

\end{document}